\documentclass[11pt]{amsart}
\usepackage{a4wide}
\usepackage{amsfonts,amsmath}
\usepackage[usenames]{color}
\usepackage[dvips]{graphicx}
\newdimen\margin   
\def\COMMENT#1{}
\let\COMMENT=\footnote

\newcommand{\eps}{\varepsilon}

\newcommand{\prob}{\mathbb{P}}
\newcommand{\ex}{\mathbb{E}}

\newcommand{\E}{{\mathcal E}}

\newcommand{\X}{{\mathcal X}}

\newcommand{\T}{{\mathcal T}}

\newcommand{\N}{{\mathcal N}}
\newcommand{\I}{{\mathcal I}}

\newcommand{\RN}{N_{t}}

\newcommand{\CN}{{\mathcal{N}}}

\renewcommand{\E}[1]{\mathbb{E}\left(#1\right)}

\newtheorem{firsttheorem}{Proposition}

\newtheorem{thm}[firsttheorem]{Theorem}
\newtheorem{lem}[firsttheorem]{Lemma}

\newtheorem{prop}[firsttheorem]{Proposition}

\numberwithin{equation}{section}
\numberwithin{firsttheorem}{section}

\begin{document}

\thispagestyle{empty}

\title{The Speed of Broadcasting in Random Networks: Density Does Not Matter}

\author{Nikolaos Fountoulakis}
\author{Anna Huber}
\author{Konstantinos Panagiotou}
\email{\{fountoul,  ahuber, kpanagio\}@mpi-inf.mpg.de}

\maketitle

\vspace{-0.6cm}
\begin{center}
\small
Max-Planck-Institute for Informatics \\
Stuhlsatzenhausweg 85, Campus E1.4 \\
Saarbr\"ucken, D-66123 Germany
\end{center} 

\begin{abstract}
Broadcasting algorithms are of fundamental importance for distributed systems engineering. In this paper we revisit the classical and well-studied \emph{push} protocol for message broadcasting. Assuming that initially only one node has some piece of information, at each stage every one of the informed nodes chooses randomly and independently one of its neighbors and passes the message to it. 

The performance of the push protocol on a fully connected network, where each node is joined by a link to every other node, is very well understood. In particular, Frieze and Grimmett proved that with probability $1-o(1)$ the push protocol completes the broadcasting of the message within $(1\pm \eps) (\log_2 n + \ln n)$ stages, where $n$ is the number of nodes of the network. However, there are no tight bounds for the broadcast time on networks that are significantly sparser than the complete graph.

In this work we consider random networks on $n$ nodes, where every edge is present with probability $p$, independently of every other edge. We show that if $p\geq {\alpha (n) \ln n \over n}$, where $\alpha(n)$ is any function that tends to infinity as $n$ grows, then the push protocol broadcasts the message within $(1\pm \eps) (\log_2 n + \ln n)$ stages with probability $1-o(1)$. In other words, in almost every network of density $d$ such that $d \ge \alpha(n) \ln n$, the push protocol broadcasts a message as fast as in a fully connected network. This is quite surprising in the sense that the time needed remains essentially \emph{unaffected} by the fact that most of the links are missing. 
\end{abstract}

\section{Introduction} 

We consider the problem of spreading information in large random networks with small average degree. Randomized broadcasting is among the most fundamental and well-studied communication primitives in distributed computing, and has also applications in several other disciplines, like e.g.\ in mathematical theories of epidemics. A particularly popular example~\cite{41841} is the maintenance of consistency in a distributed database, which is replicated at many hundreds or thousands of sites in a large, heterogeneous network. Obviously, efficient broadcasting algorithms are crucial in order to ensure that all copies of the database converge quickly and effectively to the same content.

There is an enormous amount of literature devoted to the theoretical and experimental evaluation of broadcasting algorithms on several different underlying networks. Our interest in considering \emph{random} networks is motivated, among other reasons, by P2P (peer-to-peer) systems. The idea of using random graphs appears in some ``real-life'' networks, like the popular \emph{Gnutella} network~\cite{gnutella}, or the \emph{Juxtapose} protocol~\cite{10.1109/P2P.2001.990421}, which was originally developed by Sun Microsystems. Meanwhile, a considerable amount of work by several research groups aimed at designing many diverse networks for P2P systems that resemble properties of random graphs, see e.g.~\cite{Jagannathan06queryprotocols,pru,Law03distributedconstruction}, and at developing protocols that perform efficiently on random (nearly) regular networks~\cite{1400773,eDISC,egsDISC}.

The most relevant properties of P2P networks, and more generally, of communication networks, are high expansion, connectivity, small average degree, and, (approximate) regularity of the degrees of the nodes. The random graph model considered in this paper has these properties. In particular, we investigate the classical Erd{\H{o}}s-R{\'e}nyi graph $G_{n,p}$, which is obtained by including each of the possible $\binom{n}{2}$ edges that connect any two out of $n$ labeled vertices with probability $p$, independently of all other edges.

\subsection*{The Push Model} A classical protocol in the context of randomized broadcasting, which is also the main topic of our study, is the \emph{push model}~\cite{FG85,41841}. There, initially some information is placed on one of the nodes. In each succeeding stage, every informed node passes the information to another node, that it chooses uniformly at random and independently among its neighbors. The crucial question now is: how long does it take until all nodes have received the information? There are several advantages of considering a broadcast algorithm like this: it is simple, local, and scalable, and thus independent of the network topology. Moreover, it is highly robust against network and link failures, which makes it highly reliable.

In the case where the underlying network is the complete graph on~$n$ vertices, Frieze and Grimmett~\cite{FG85} proved that with high probability\footnote{with probability tending to 1 when~$n\to\infty$} (w.h.p.) the push protocol completes the broadcasting of the message within~$(1\pm \eps) (\log_2 n + \ln n)$ stages. In other words, if a node can ``talk'' to \emph{any other} node in the network, then the broadcast time will be almost surely very close to~$\log_2 n + \ln n$. This bound was later improved by Pittel~\cite{Pi87} to~$\log_2 n + \ln n + \alpha(n)$, where~$\alpha(n)$ is any function that tends to~$\infty$ when~$n\to\infty$. Feige et \ al.\ considered in~\cite{FPRU90} networks that are different from the complete graph. Among other results, they showed that if the underlying network is a random graph~$G_{n,p}$, where~$p\ge \frac{(1+\eps)\ln n}{n}$, then the message will arrive at all nodes with high probability within~$\Theta(\ln n)$ stages. Moreover, they also showed that the protocol is efficient on hypercubes, and derived bounds that hold for arbitrary graphs. Els\"asser and Sauerwald determined in \cite{ES07} similar bounds for several classes of Cayley graphs, thus generalizing upon~\cite{FPRU90}.

\subsection*{Our contribution}
Let $G = (V,E)$ be a graph on $n$ vertices, where we will assume that  $V=\{1,\ldots, n\}$. We define $T(G)$ as the number of stages needed by the push protocol until all vertices have been informed, if the information 
is initially placed on node 1. In the remainder of the paper, we will be using the terms ``node" and ``vertex" without distinction. 
Note that \emph{regardless} of the underlying network topology~$T(G)\geq \log_2 n$, as the number of informed vertices can at most double in each round. Consequently, all the results mentioned above state that the push model is, up to multiplicative constants, an asymptotically optimal protocol for disseminating information.

However, it is not at all well-understood \emph{how much} the structure of the underlying network affects the performance of the push model. Although, for example, we know from the results in~\cite{FPRU90} that on a random graph $G_{n,p}$ the protocol requires with high probability at most $C\ln n$ rounds, for some $C > 0$, we have \emph{a priori} no bounds than quantify how slower (or faster?) the protocol is compared to the case where the network is the complete graph. In particular, it is not clear in which way the average degree of the underlying graph influences the speed of the protocol. Our main result states that the number of stages is essentially \emph{unaffected} by the density of the underlying graph, thus confirming the robustness and the efficiency of the push model:
\begin{thm}
\label{thm:main}
Let $0 < \alpha(n)\leq \ln^{1/9} n$ be any function with the property $\lim_{n\to\infty} \alpha(n) = \infty$. Let $p \ge \frac{\alpha(n) \ln n}{n}$. Then w.h.p.
\[
	|T(G_{n,p}) - (\log_2 n + \ln n)| < {\alpha (n)^{-1/7} \ln n}.
\]
\end{thm}
In other words, if the average degree of $G_{n,p}$ is slightly larger than $\ln n$, then the broadcast time of the push model essentially \emph{coincides} with the broadcast time on the complete graph, which was shown in~\cite{FG85} to be very close to $\log_2 n + \ln n$. Consequently,
the number of stages needed is not influenced by the fact that most of the links are missing. 

To avoid any confusion, we want to note that in Theorem~\ref{thm:main} the term ``w.h.p'' refers to two independent probability spaces: first, the space from which we sample the underlying network, and second, the space of the random choices performed by the nodes.

\subsection*{Proof Ideas \& Techniques}
Before we proceed with a detailed exposition of our proof, let us mention a few words about the general strategy. Theorem~\ref{thm:main} is proved by
bounding for each stage performed by the push model simultaneously from above and and from below the number of informed nodes. In particular, we show
that in the first $(1 - o(1))\log_2 n$ stages, the number of informed nodes nearly \emph{doubles} in each stage. As a result, we are able to show that
after nearly $\log_2 n$ rounds there will be $\eps n$ informed nodes in total. Then things evolve very fast: only after a small number of stages, the
number of nodes having the information will be already roughly $(1 - \eps)n$. After that, we show that additionally approximately $\ln n$ stages are
necessary and sufficient to spread the information to everybody.

The analysis of the last stages is particularly challenging from a technical point of view, as the number of informed nodes increases only slowly towards the end of the process. In such cases, it is typically difficult to control the deviations of several involved random variables from their expectations. To this end, we exploit a modern and powerful tool from probability theory called \emph{Talagrand's inequality}, which -- to our knowledge -- has not been applied in the context of distributed computing problems. We believe that it could be widely applicable to the analysis of existing or future randomized protocols with several different degrees of dependency.

\subsection*{Outline}
Section \ref{sec:preliminaries} introduces the main tools from probability theory that we will use, and in particular Talagrand's inequality. In
Section \ref{sec:propGnp} we collect and prove the basic properties of $G_{n,p}$ that will be important in the proof of Theorem~\ref{thm:main}, and
introduce some necessary notation that will be used throughout. Finally, Section~\ref{sec:broadcasting} contains the ``core'' of the proofs, where
the general strategy given above is converted to  a rigorous argument.


\section{Preliminaries} 
\label{sec:preliminaries}
A basic tool that we will use in the following proofs is the Chernoff bound. 
This provides exponentially small bounds for the probability that a binomially 
distributed random variable deviates significantly from its expected value. 
\begin{thm}[Chernoff Bounds]
\label{thm:chernoff}
Let $X$ be a binomially distributed random variable and let $x> 0$. Then 
$$ \prob (|X-\ex (X)|>x)\leq 2 \exp \left(- {x^2 \over 2(\ex (X) + x/3)}\right).$$
\end{thm}
A more general tool that we shall apply several times is the inequality by Azuma and Hoeffding. Intuitively, it provides strong bounds for the probability that a function defined on a set of independent random variables deviates significantly from its expectation, when the value of the function is not affected much by small changes in each one of its arguments.
\begin{thm}[Azuma-Hoeffding's Inequality]
\label{thm:azuma}
Let $Z_1,\ldots, Z_N$ be independent random variables taking values in the sets $\Lambda_1, \ldots , \Lambda_N$ respectively. Let $\Lambda = \Lambda_1 \times \cdots \times \Lambda_N$. Let $f:\Lambda \rightarrow \mathbb{R}$ be a function and set $X = f(Z_1,\ldots, Z_N)$. Assume that there are quantities $c_k$, $k=1,\ldots , N$ satisfying the 
following: 
\begin{enumerate}
\item [a.] If $z,z' \in \Lambda $ differ only in the $k$th coordinate, then $|f(z) - f(z')|\leq c_k$.
\end{enumerate} 
Then, for every $x\geq 0$ we have that
\begin{equation} \label{eq:azuma}
\prob (|X-\mathbb{E}(X)| \geq x) \leq 2 \exp \left(- {x^2 \over 2\sum_{i=1}^Nc_i^2} \right).
\end{equation}
\end{thm}
Note that the above inequality gives meaningful bounds only if the expectation of $X$ is much larger than $(\sum_{i=1}^N c_i^2)^{1/2}$. This condition is unfortunately not always given in our intended applications. In such cases, we will use an estimate given by Talagrand (see the following theorem), which gives a much stronger tail bound, provided that an additional assumption is satisfied. Intuitively, the statement claims that if the value of $X$ is ``witnessed" by only a ``small'' number of its arguments, then $X$ is sharply concentrated. However, there is a small caveat: the concentration is not guaranteed to be around the expectation, but instead around the \emph{median} of $X$. (Recall that the median is a number $m$ such that $\prob (X<m) \leq {1\over 2}$ and $\prob (X > m)\leq {1\over 2}$.) 
As we shall see below,  this is not a significant problem as in general the median is very close the expected value. 
\begin{thm}[Talagrand's Inequality]
\label{thm:talagrand}
Suppose that the preconditions of Theorem~\ref{thm:azuma} are satisfied. Additionally, assume that there is an increasing function $\psi$ satisfying the 
following: 
\begin{enumerate}
\item [b.] Let $z \in \Lambda$ and $r \in \mathbb{R}$ such that $f(z) \geq r$. Then there exists a set $J \subseteq \{1,\ldots, N \}$ with $\sum_{i \in J}c_i^2 \leq \psi (r)$, such that for all $y \in \Lambda$ with $y_i = z_i$ when $i \in J$, we have $f(y) \geq r$.  
\end{enumerate} 
Then, if $m$ is the median of $X$, for every $x\geq 0$ we have 
\begin{equation} \label{eq:Talagrand}
\prob (|X-m|\geq x) \leq 4 \exp \left(- {x^2 \over 4\psi (m + x)} \right).
\end{equation}
\end{thm}
The next statement gives a sufficient condition that ensures that the median is very close to the expected value.
\begin{prop}[Example 2.33 in~\cite{JLR}]
Let $X$ be a random variable that satisfies the preconditions of Theorem~\ref{thm:talagrand} with $\psi(r) \le \lceil r\rceil$. Then
\begin{eqnarray} \label{eq:MedMean}
|m - \ex (X)| = O \big(\sqrt{\ex (X)} \big).
\end{eqnarray}
\end{prop}
The presentation of the above inequalities is as in~\cite{JLR}, where also many applications are presented.


\section{Properties of $G_{n,p}$} 
\label{sec:propGnp}

For any graph $G$ with vertex set $V$ let $\Gamma_G(v)$ be the set of neighbors of $v$ in $G$. Moreover, for $S, S' \subseteq V$ we will denote by $e_G(S,S')$ the number of edges with precisely one endpoint in each of $S, S'$. Finally, for two real numbers $a,b$ we will write $a \pm b$ for the interval of reals $(a-b, a+b)$, and with slight abuse of notation we will write $X = a\pm b$ to denote $X \in a\pm b$.

Let $\alpha(n) > 0$ be any function with $\lim_{n\to\infty} \alpha(n) = \infty$ and let~$p \ge \frac{\alpha(n) \ln n}{n}$.
In this section we collect a few properties of $G_{n,p}$ that we will use in the proof of Theorem~\ref{thm:main}. 

Note that for any $S \subset V$, the expected number of neighbors of any $v \in V\setminus S$ in $S$ is $p|S|$. The next lemma says that \emph{for all} large enough $S$ almost all vertices have roughly the right degree in $S$.
\begin{lem}
\label{lem:bigExceptional}
 The random graph~$G_{n,p}$ has w.h.p.\ for any~$\alpha(n)^{-1/2}\le \eps\le 1$ the following property. 
For any subset~$S$ of its vertices satisfying~$|S| \ge \frac{n}{\alpha(n)}$, there is a set~$X_S \subset V\setminus S$ that contains at most~$\frac{8n}{\ln n}$ vertices such that
\[
	\forall v\in (V\setminus S)\setminus X_S : |\Gamma_{G_{n,p}}(v) \cap S| =  (1\pm \eps) p |S|.
\]
\end{lem}
\begin{proof}
Let~$S$ be any fixed subset of the vertices such that~$|S| \ge \frac{n}{\alpha(n)}$. We call a vertex~$v \in V\setminus S$ \emph{violating} with respect to~$S$, if the number of its neighbors in~$S$ is~$> (1+\eps)p|S|$ or~$<(1-\eps)p|S|$. Assume there exist at least~$t := \frac{8n}{\ln n}$ vertices that are violating, and denote by~$X_S$ the set consisting of those vertices.

Note that the expected number of neighbors in~$S$ of a vertex is~$p|S|$. By applying the Chernoff bounds, we obtain that the probability that a vertex is violating is for large~$n$ at most~$e^{-\eps^2 p|S|/4}$. Moreover, the events that two distinct vertices are violating are independent, which implies that the probability that there are~$t$ violating vertices is bounded from above by~$e^{-\eps^2 p |S|/4 \cdot t}$. Hence, as there are~$\binom{n}{|S|} \le n^{|S|} = e^{|S|\ln n}$ ways to choose~$S$, the probability that there is a set such that there are~$t$ violating vertices with respect to it as at most
\[
	\exp\left\{|S|\ln n - \frac{\eps^2 p |S|}4 \cdot t\right\}
	\le \exp\left\{|S|\left(\ln n - \frac{\eps^2 p}4 \cdot \frac{8n}{\ln n}\right)\right\}.
\]
This, combined with the bound~$p \ge \frac{\alpha(n)\ln n}{n}$, can be estimated with plenty of room to spare from above by at most~$e^{-|S|\ln n}$. The proof is completed by summing this expression up for all~$|S| \ge \frac{n}{\alpha(n)}$.
\end{proof}

The next statement considers a similar setting as before, with the difference that now $S$ might be very small. Here we show that the number of vertices that have many neighbors in~$S$ is only $o(|S|)$.
\begin{lem}
\label{lem:smallExceptional}
For any $\eps \ge \alpha(n)^{-1/2}$, the random graph $G_{n,p}$ has w.h.p.\ the following property.
For any subset $S$ of its vertices such that $|S| \le \frac{n}{\alpha(n)}$ there is a set $X_S$ containing at most~$|S|\eps^{-1}\alpha(n)^{-1}$ vertices, such that
\[
	\forall v\in (V\setminus S)\setminus X_S : |\Gamma_{G_{n,p}}(v) \cap S| \le \eps p n.
\]
\end{lem}
\begin{proof}
The proof is similar to the proof of Lemma~\ref{lem:bigExceptional}, except that here we have to deal with small sets $S$. 
We give the whole proof for the sake of completeness. We assume that $|S| \geq \eps p n$, for otherwise the statement holds trivially.

Let~$S$ be any fixed subset of the vertices such that~$|S| \le \frac{n}{\alpha(n)}$. We call a vertex~$v \in V\setminus S$ \emph{violating} with respect to~$S$, if the number of its neighbors in~$S$ is~$>\eps p n$. Assume there exist at least~$t := \frac{|S|}{\eps \alpha(n)}$ vertices that are violating, and denote by~$X_S$ the set consisting of those vertices.

The expected number of neighbors in $S$ of a vertex $v \in V\setminus S$ is $p|S| = o(\eps pn)$. A straightforward application of the Chernoff bounds then implies that the probability that a vertex is violating is for large $n$ at most $e^{-\eps p n}$. Hence, as the events that distinct vertices are violating with respect to $S$ are independent, the probability that there are $t$ such vertices is at most $e^{-\eps pn \cdot t}$.

Note that the number of ways to choose $S$ is $\binom{n}{|S|} \le (\frac{en}{|S|})^{|S|}$. In conclusion, the probability that there is an $S$ with $t$ violating vertices is at most
\[
	\left(\frac{e}{|S|}\right)^{|S|}\exp\left\{|S|\ln n - \eps p n \cdot t\right\}
	\le \left(\frac{e}{|S|}\right)^{|S|} \exp\left\{|S|\left(\ln n - p\alpha(n)^{-1} \right)\right\}
	\le \left(\frac{e}{|S|}\right)^{|S|}.
\]
The proof then completes by summing this expression up for all~$\eps p n\le |S| \le \frac{n}{\alpha(n)}$.
\end{proof}
Finally, we need the following statement about the distribution of the edges in~$G_{n,p}$. The proof is a straightforward application of Chernoff's
bounds, and quite standard in the classical random graph theory. We include a short proof for completeness.
\begin{lem}
\label{lem:degsAndEdges}
The following holds w.h.p.
\[
	\forall S \subseteq V: e_{G_{n,p}}(S, V\setminus S) = |S|(n - |S|)p \left(1 \pm \sqrt8\alpha(n)^{-1/2}\right).
\]
\end{lem}
\begin{proof}
It is sufficient to show the statements for $S$ such that $|S| \le n/2$. For any fixed such $S$, the quantity $e_{G_{n,p}}(S, V\setminus S)$ is
binomially distributed with expectation $|S|(n-|S|)p$. Call $S$ \emph{bad}, if $e_{G_{n,p}}(S, V\setminus S)$ deviates from its expected value by more than 
$\sqrt{4|S|^2(n - |S|)p\ln n}$. Note that
$$ {\sqrt{4|S|^2(n - |S|)p\ln n}\over |S|(n - |S|)p} = \sqrt{4\ln n \over np(1-|S/n|)} \leq \sqrt{8 \ln n \over np}\leq \sqrt{8 \over \alpha(n)}.$$
By applying the Chernoff bounds we obtain that the probability that $S$ is bad is with plenty of room to spare for large~$n$ at most
\[
	\exp\left\{-\frac{4|S|^2(n-|S|)p \cdot \ln n}{3|S|(n-|S|)p}\right\}
	= \exp\left\{-\frac43|S| \ln n \right\}.
\]
Then, as the number of ways to choose $S$ is at most~$n^{|S|}$, we infer by summing over all $1 \le |S| \le n/2$ that w.h.p.\ there is no bad set $S$ in $G_{n,p}$. The proof completes readily by using that $n - |S| \ge n/2$ and the lower bound on $p$.
\end{proof}
Note that in the special case that $|S| = 1$ in the above lemma, i.e., $S$ contains just a single vertex $v$, we obtain that
\begin{equation*}
	|\Gamma_{G_{n,p}}(v)| = e_{G_{n,p}}(S, V\setminus S) = (1\pm 3\alpha(n)^{-1/2})pn.
\end{equation*}
This fact will become very handy later and we will use it without further reference.

\section{Broadcasting on Random Graphs} 
\label{sec:broadcasting}

Let $G$ be any graph with vertex set $V$ and let $p\ge {\alpha (n) \ln n \over n}$, where $\alpha(n) \le \ln^{1/9}n$ is any positive function such that $\lim_{n\to\infty}\alpha(n) = \infty$. Fix 
\[
	\eps := \alpha(n)^{-1/2}.
\]
We say that $G$ is $p$-\emph{typical} if it satisfies the following three conditions:
\begin{itemize}
	\item[(I)] For any $S \subseteq V$ such that $|S| \ge \frac{n}{\alpha(n)}$ there is a $X_S \subset V \setminus S$ such that $|X_S| \le \frac{8n}{\ln n}$ and
		\[
			\forall v\in (V\setminus S)\setminus X_S : |\Gamma_{G}(v) \cap S| =  (1\pm \eps) p |S|.
		\]
	\item[(II)] For any $S \subseteq V$ such that $|S| \le \frac{n}{\alpha(n)}$ there is a $X_S \subset V \setminus S$ such that $|X_S| \le \frac{|S|}{\eps \alpha(n)}$ and
		\[
			\forall v\in (V\setminus S)\setminus X_S : |\Gamma_{G}(v) \cap S| \le \eps p n.
		\]
	\item[(III)] For all $S \subseteq V$
		\[
			e_{G}(S, V\setminus S) = |S|(n - |S|)p \left(1 \pm \sqrt8\eps\right).
		\]
\end{itemize}
We will denote by $\T_n(p)$ the set of $p$-typical graphs on $V$. Note that Lemmas~\ref{lem:bigExceptional}-\ref{lem:degsAndEdges} guarantee that $G_{n,p}$ is w.h.p.\ $p$-typical. Hence, we shall restrict our attention only to graphs in $\T_n(p)$.

~\\
Let us denote by $T_1(G)$ the first point in time where at least $\eps n$ vertices are informed and $T_2(G)$ the first point in time where at least $(1-\eps)n$ vertices are informed.  Our aim is to give bounds on $T(G)$ by bounding $T_1(G)$, $T_2 (G) -  T_1(G)$ and $T(G) - T_2(G)$ uniformly for every $G \in \T_n(p)$. The following three lemmas do so. In the proofs we will several times assume that $n$ is sufficiently large so that the claimed inequalities hold, without explicitly mentioning that.
\begin{lem}
\label{lem:doubling}
Uniformly for $G\in \T_n(p)$, with probability $1 - o(1)$ it holds that $$|T_1(G) - \log_2 n| \le 9\sqrt{\eps}\log_2 n.$$
\end{lem}
\begin{proof}
First of all, note that always $T_1(G) \ge \log_2(\eps n)$, as the number of informed nodes at most doubles in each stage. Hence, we restrict our attention to the proof of the upper bound for $T_1(G)$.

Let $\I_t$ be the random set of informed vertices after $t$ stages, and set $I_t := |\I_t|$. Note that our definitions imply that $\I_0 = \{1\}$. We will show that
\begin{equation}
\label{eq:almost2x}
	\prob\big(I_{t+1} \ge (2 - 7\sqrt{\eps})I_t ~|~ I_t < \eps n\big)
	\ge
	1 -
	\begin{cases}
	o\left((\ln n)^{-1}\right) ,& I_t \ge \ln^{1/4} n\\
	\ln^{-1/2}n,& \text{otherwise}
	\end{cases}.
\end{equation}
The proof of the lemma then completes by a repeated application of the above inequality. In particular, either there is a $t < (1 + 8\sqrt{\eps})\log_2 n$ such that $I_t \ge \eps n$, in which case there is nothing to show, or, with probability $1 - o(1)$,
\[
	I_{\lceil (1 + 8\sqrt{\eps})\log_2 n \rceil} \ge (2-7\sqrt{\eps})^{(1 + 8\sqrt{\eps})\log_2 n} \ge \eps n.
\]
So we showed that
\[
T_1(G) \leq \lceil (1 + 8\sqrt{\eps})\log_2 n \rceil \leq (1 + 9\sqrt{\eps})\log_2 n
\]
In the remainder we prove~\eqref{eq:almost2x}. For every vertex $v\in I_t$ we define an indicator random variable $N_v$ that equals 1 if $v$ informs a vertex in $V \setminus \I_t$. Moreover, for every pair of distinct vertices $v,v' \in \I_t$ let $C_{v,v'}$ be the indicator variable that is equal to 1 if $v$ and $v'$ inform the same vertex in $V \setminus \I_t$. Finally, denote by $\CN_t$ the random set of vertices in $V\setminus \I_t$ that will be informed in stage $t+1$ by the vertices in $\I_t$. By simple inclusion-exclusion we obtain that
\begin{equation*}
	|\CN_t|
	~\ge~
	\sum_{v\in \I_t} N_v - \sum_{v,v'\in \I_t, v\neq v'} C_{v,v'}.
\end{equation*}
Note that
\begin{equation}
\label{eq:auxNC}
	\E{N_v} = \frac{|\Gamma_G(v) \cap (V \setminus \I_t)|}{|\Gamma_G(v)|}
	~\text{ and }~
	\E{C_{v,v'}} = \frac{|\Gamma(v) \cap \Gamma(v') \cap (V \setminus \I_t)|}{|\Gamma(v)||\Gamma(v')|}.
\end{equation}
We shall now show that $|\CN_t| \ge (1 - 7\sqrt{\eps})I_t$ holds with the desired probability, which will complete the proof of~\eqref{eq:almost2x}. To achieve this we shall argue differently in the two cases $I_t \ge \ln^{1/4} n$ and $I_t < \ln^{1/4} n$. Before we proceed, let us make two auxiliary preparations. Note that by property (III) of $G$ we obtain for sufficiently large $n$ that
\[
	\forall v \in V: |\Gamma_G(v)| = (1\pm3\eps)pn.
\]
This, together with \eqref{eq:auxNC} implies with a simple double counting argument that
\begin{equation*}
\label{eq:sum}
	\sum_{\substack{v,v' \in \I_t \\ v \not = v'}} \E{C_{v,v'}}
	=
	\sum_{\substack{v,v' \in \I_t \\ v \not = v'}} \frac{|\Gamma_G(v) \cap \Gamma_G(v') \cap \left( V \setminus \I_t \right)|}{(1\pm7\eps)(pn)^2}
	=  
	\frac{(1\pm8\eps)}{(pn)^2} \cdot \sum_{u \in V \setminus \I_t} \binom{|\Gamma_G(u) \cap \I_t|}2. 
\end{equation*}
We will use these facts in the remainder without further reference.

~\\
First, suppose that $I_t < \ln^{1/4} n$. Note that for each vertex $v \in \I_t$ at least $|\Gamma_G(v)|- \ln^{1/4} n$ of the edges that are adjacent to it are directed to vertices outside $\I_t$. This implies that
$$
	\prob (N_v =1)
	=
	\frac{|\Gamma_G(v) \cap (V \setminus \I_t)|}{|\Gamma_G(v)|}
	\ge
	1 - { \ln^{1/4} n \over (1-4\eps)pn}
	\ge 1 - {1\over 2}~\ln^{-3/4} n.
$$
Therefore, with probability at least $1-{1\over 2}~\ln^{-1/2}n$, all vertices in $\I_t$ inform a vertex that lies outside $\I_t$, i.e., $\sum_{v\in \I_t} N_v = I_t$. However, there is still the possibility that two vertices in $\I_t$ inform the same vertex, thus creating a conflict. The probability that such a conflict occurs is for large $n$ smaller than
\begin{equation*}
\begin{split}
\sum_{v,v'\in \I_t, v\neq v'} \E{C_{v,v'}}
\le
\frac{2}{(pn)^2} \cdot \sum_{u \in V \setminus \I_t} \binom{|\Gamma_G(u) \cap \I_t|}2.
\end{split}
\end{equation*} 
Note that $0 \le |\Gamma_G(u) \cap \I_t| \le I_t$. Moreover, property (III) in the definition of $\T_n(p)$ implies that
\[
	\sum_{u\in V\setminus \I_t} |\Gamma_G(u) \cap \I_t|
	= e_G(\I_t, V\setminus \I_t)
	\le 2I_tpn.
\]
Under these conditions, as the sum of the binomial coefficient above is a convex function, it is bounded from above when we set $|\Gamma_G(u) \cap \I_t| = I_t$ for $2pn$ choices of $u$, and $|\Gamma_G(u) \cap \I_t| = 0$ otherwise. Hence, we obtain for large $n$ that
\[
	\sum_{v,v'\in \I_t, v\neq v'} \E{C_{v,v'}}
	\le \frac{2}{(pn)^2}\cdot 2pn \cdot I_t^2
	\le {1\over 2}~\ln^{-1/2} n.
\]
So, with probability at least $1 - {1\over 2}~\ln^{-1/2}n - {1\over 2}~\ln^{-1/2} n \ge 1 - \ln^{-1/2}n$ we have that $|\N_t| = I_t$, which completes the proof for the case $I_t < \ln^{1/4} n$.

~\\
Finally, we consider the case $I_t \ge \ln^{1/4} n$. We will first give tight bounds on the expectation of $|\CN_t|$, and then apply the Azuma-Hoeffding inequality to show that $|\CN_t|$ is sufficiently sharply concentrated around $\E{|\CN_t|}$. By using~\eqref{eq:auxNC} we obtain with plenty of room to spare for large~$n$ that
\begin{equation}
\label{eq:sumNv}
	\E{\sum_{v\in \I_t} N_v}
	= \sum_{v\in \I_t} \frac{{|\Gamma_G(v) \cap (V \setminus \I_t)|}}{(1\pm3\eps)pn}
	= \frac{(1\pm4\eps)e_G(\I_t, V\setminus \I_t)}{pn}
	\stackrel{(III)}{=} (1\pm 8\eps)I_t. 
\end{equation}
Moreover, recall that
\begin{equation}
\label{eq:TMPbinomials}
	\sum_{{v,v' \in \I_t, v \not = v'}} \E{C_{v,v'}}
	=  
	\frac{(1\pm8\eps)}{(pn)^2} \cdot \sum_{u \in V \setminus \I_t} \binom{|\Gamma_G(u) \cap \I_t|}2. 
\end{equation}
We are going to estimate the last sum from above as follows. As $G\in\T_n(p)$ we may infer the following.
\begin{itemize}
	\item If $I_t \le \frac{n}{\alpha(n)}$, then, by (II), there is $\X \subset V\setminus \I_t$ such that $|\X| \le \sqrt{\eps}I_t$ and
		\[
			\forall v \in (V \setminus \I_t) \setminus \X : |\Gamma_G(v) \cap \I_t| \le \eps pn.
		\]
	\item If $\frac{n}{\alpha(n)}\le I_t \le \eps n = \frac{n}{\alpha(n)^{1/2}}$, then, by (I), there is $\X \subset V\setminus \I_t$ such that $|\X| \le \frac{8 n}{\ln n}$ and
		\[
			\forall v \in (V \setminus \I_t) \setminus \X : |\Gamma_G(v) \cap \I_t| \le (1 + \eps) pI_t \le 2\eps pn.
		\]
\end{itemize}
So, in both cases we have for all $v \in (V \setminus \I_t) \setminus \X$ that $|\Gamma_G(v) \cap \I_t| \le 2\eps pn$, and $|\X| \le \sqrt{\eps}I_t$. Moreover, by exploiting property (III) of $G$ we obtain that for all $v\in\X$ it holds $|\Gamma_G(v) \cap \I_t| \le 2pn$. Using this, we can bound from above the sum in~\eqref{eq:TMPbinomials} by splitting it into two parts as follows:
\begin{equation*} 
\begin{split}
&\sum_{u \in V \setminus \I_t} |\Gamma_G(u) \cap \I_t|^2 = \sum_{u \in (V \setminus \I_t)\setminus \X} |\Gamma_G(u) \cap \I_t|^2 + 
\sum_{u \in \X} |\Gamma_G(u) \cap \I_t|^2 \\
& \leq \sum_{u \in (V \setminus \I_t)\setminus \X}  |\Gamma_G(u) \cap \I_t|^2 + |\X|~(2pn)^2 \leq
\sum_{u \in (V \setminus \I_t)\setminus \X}  |\Gamma_G(u) \cap \I_t|^2 + \sqrt{\eps}I_t ~(2pn)^2.
\end{split}
\end{equation*}
Note that $0 \le |\Gamma_G(u) \cap \I_t|  \leq 2 \eps p n$ for every $u \in (V \setminus \I_t)\setminus \X$. Moreover, it is easily seen that $\sum_{u \in V \setminus \I_t} |\Gamma_G(u) \cap \I_t| = e_G(\I_t, V \setminus \I_t)$. By the convexity of $x^2$, the sum in the expression above is bounded from above if we choose $|\Gamma_G(u) \cap \I_t| = 2\eps pn$ for $e_G(\I_t, V \setminus \I_t)/(2\eps p n)$ different $u$'s, and $|\Gamma_G(u) \cap \I_t| = 0$ otherwise. We obtain
\begin{equation*} 
\begin{split}
\sum_{u \in V \setminus \I_t} |\Gamma_G(u) \cap \I_t|^2
\leq {I_t (n - I_t)p (2 \eps p n)^2 \over \eps pn} + \sqrt{\eps} I_t ~(2pn)^2
\leq 
\frac{9}{2}\sqrt{\eps}p^2n^2I_t. 
\end{split}
\end{equation*}
By plugging this into~\eqref{eq:TMPbinomials} we obtain that $	\sum_{{v,v' \in \I_t, \ v \not = v'}} \E{C_{v,v'}} \le 5\sqrt{\eps}I_t$. Finally, combined with~\eqref{eq:sumNv} this gives with lots of room to spare that.
\[
	\E{|\N_t|} \ge (1 - 6\sqrt{\eps})I_t.
\]
To complete the proof we will bound the probability that $|\N_t| < I_t (1 - 7\sqrt{\eps})$ by using Azuma-Hoeffding's inequality. Note that $|\N_t|$ can
change by at most 1, if we modify one of the choices made by some vertex in $\I_t$. So, by applying Theorem~\ref{thm:azuma} with $c_i = 1$ and 
$N = I_t$ we obtain
\begin{equation*} 
\begin{split}
\prob & \left(|\N_t| < I_t(1 - 7 \sqrt{\eps}) \right)
\le
\prob \left( |\N_t| < \ex \left( |\N_t|\right) - \sqrt{\eps}I_t \right)
\le e^{-\frac{\eps}{2}\ln^{1/4}n},
\end{split}
\end{equation*}
thus concluding the proof for the case $I_t \ge \ln^{1/4}n$.
\end{proof}
In the next lemma we will consider the ``intermediate'' phase of the push model between~$T_1 (G)$ and~$T_2 (G)$ for~$G \in \T_n(p)$. Our general strategy is to bound the number~$N_t$ of vertices that get informed in the current stage~$t$ from below. For this, we first estimate~$\ex (\RN)$ and then we use concentration inequalities (Theorem~\ref{thm:azuma}) to show that with sufficiently high probability~$\RN$ is very close to~$\ex (\RN)$.
\begin{lem}
\label{lem:slowDown}
Uniformly for all $G \in \T_n (p)$,  with probability $1-o(1)$ it holds that 
$$T_2(G) - T_1(G) \leq  {9\eps^{-1}\ln \eps^{-1}},$$
and there are at least $\eps n /(2e)$ uninformed vertices at $T_2(G)$. 
\end{lem}
\begin{proof} 
Let $\I_t$ be the random set of informed vertices after $t$ stages, and set $I_t := |\I_t|$. We will show that for $T_1(G) \leq t < T_2(G)$
\begin{equation} \label{eq:RecursioToProve}
I_{t+1} \geq I_t \left( 1 + {\eps \over 4} \right),
\end{equation}
with probability at least $1- e^{-\eps^6 n /8}$. 
Let us abbreviate $b = 8 \eps^{-1}\ln \eps^{-1}$. To see that this is sufficient, note that if ``$T_2 (G) - T_1 (G) \leq  b$", then there is 
nothing to prove. On the other hand, if ``$T_2 (G) - T_1 (G) > b$", then with (conditional) probability at least $(1- e^{-\eps^6 n /8})^{b} = 1- o(1)$, 
for $\lceil b\rceil$ consecutive steps after $T_1(G)$ the recursion (\ref{eq:RecursioToProve}) holds.
In turn, this implies with $1+x > e^{x/2}$, which is valid for small enough $x>0$, that
\begin{equation*}
I_{T_1(G) +\lceil b\rceil} 
\geq I_{T_1(G)} \cdot \left( 1 + {\eps \over 4} \right)^{\lceil b \rceil} 
> \eps n e^{b\eps/8} > (1-\eps)n.
\end{equation*}
Therefore $I_{T_1(G) + \lceil b\rceil} > n(1-\eps)$, from which we obtain with plenty of room to spare that, say, $T_2 (G) - T_1 (G) \leq b+1\leq 
{9 \eps^{-1}\ln \eps^{-1}}$.

Now we turn to the proof of~\ref{eq:RecursioToProve}.
Let $t$ be such that $T_1 (G) \leq t < T_2(G)$, and denote by $\N_t$ the set of vertices in $V\setminus \I_t$ that will be informed by the vertices in
$\I_t$ in stage $t+1$. Moreover, write $N_t := |\N_t|$. We will show that $N_t$ is not much smaller than its expected value. But first let us
calculate $\E{N_t}$. The definition of the push model implies that the probability that any $v \in V\setminus \I_t$ \emph{does not} belong to $\N_t$
is precisely
\[
	\prod_{u \in \Gamma_G(v) \cap \I_t} \left(1 - \frac1{|\Gamma_G(u)|}\right).
\]
Next we make use of property (I) in the definition of $\T_n(p)$: All vertices
in $V \setminus \I_t$, apart from an exceptional set $\X = \X_t \subset V\setminus \I_t$ that contains at most $8n/\ln n$ vertices, have $(1 \pm \eps)pI_t $
neighbors in $\I_t$. 
Using this and the above fact we may write
\begin{equation} \label{eq:ExpNewlyInf}
\begin{split}
\ex \left( N_t \right) & = 
\sum_{v \in (V \setminus \I_t)\setminus \X} \left(1 - \prod_{u \in \Gamma_G(v) \cap \I_t} \left(1 - {1\over |\Gamma_G(u)|} \right) \right) 
\pm |\X|.
\end{split}
\end{equation} 
Next we derive tight bounds for the product above. Firstly, observe that property (III) implies that for all $u \in \I_t$ 
\begin{equation} \label{eq:degsLoose}
|\Gamma_G (u)| = np \left( 1 \pm 3 \eps \right). 
\end{equation}
Also, the definition of $\X$ implies for $v \in  (V \setminus \I_t)\setminus \X$
\begin{equation} \label{eq:degsinS}
|\Gamma_G (v) \cap \I_t| = \left( 1 \pm \eps \right)pI_t.
\end{equation}
Recall that for $x > 0$ small enough we have $e^{-x - x^2}\leq 1-x\leq e^{-x}$. 
So the bounds in (\ref{eq:degsLoose}) and (\ref{eq:degsinS}) imply that for $v \in (V \setminus \I_t)\setminus \X$ we have
\begin{equation*} \label{eq:prod}
\prod_{u \in \Gamma_G (v) \cap \I_t} \left(1 - {1\over |\Gamma _G(u)|} \right) = 
\exp \left( - {I_t\over n}\left(1 \pm 5\eps \right)\right) \left( 1+ O \left({1\over np} \right) \right).
\end{equation*}
As $|(V \setminus \I_t) \setminus \X| = (n-I_t)(1\pm\eps)$,
by substituting the above estimate into (\ref{eq:ExpNewlyInf}) we obtain
\begin{equation} \label{eq:ExpNewlyInf1}
\begin{split}
\ex \left( N_t \right) & = n\left(1 - {I_t \over n} \right) \left(1 - e^{-{I_t\over n}} \right) \left(1 + O (\eps)\right).
\end{split}
\end{equation}
We will bound the probability that $|\RN - \ex \left( \RN \right)| > \eps\ex \left( \RN \right)$ using the Azuma-Hoeffding 
inequality. Firstly, note that as $\eps <{I_t\over n} \leq 1-\eps$, we have 
$\ex \left( \RN \right)  \geq \eps^2 n /2$, for $n$ sufficiently large. 
Moreover, if we change only one of the random choices of the vertices in $\I_t$, then $\RN$ changes by at most 1. Thus, a simple application of Theorem~\ref{thm:azuma} with $c_k = 1$ and $N = I_t$ yields
\begin{equation*} 
\prob 
\left(|\RN - \ex \left( \RN \right)| > \eps  \ex \left( \RN \right) \right) \leq 2\exp \left( - {\eps^2\ex^2 \left( \RN \right) \over 2I_t}\right)
\leq 2\exp \left(- {\eps^6 n \over 8} \right). 
\end{equation*}
So, for $n$ sufficiently large, with plenty of room to spare we obtain that 
\begin{equation} \label{eq:ILow} 
\RN  = n\left(1 - {I_t \over n} \right) \left(1 - e^{-{I_t\over n}} \right) \left(1 \pm \sqrt{\eps}\right)
\end{equation}
with probability $2\exp \left(- {\eps^6 n \over 8} \right)$.
This identity enables us to write a recursive formula concerning the evolution of the number of informed vertices. Recall that for all $0< x <1$, we have $1-e^{-x} \geq x/2$. (\ref{eq:ILow}) implies that 
\begin{equation} \label{eq:Recursion} 
\begin{split}
I_{t+1} &\ge I_t +   n\left(1 - {I_t \over n} \right) {I_t\over 2n}\left(1 - \sqrt{\eps}\right)  
= I_t \left(1 + {1\over 2}~\left(1 - {I_t \over n} \right) \left(1 - \sqrt{\eps}\right)\right). 
\end{split}
\end{equation}
Since $I_t \leq (1-\eps) n$, it follows that for $n$ large enough
$$ {1\over 2}~\left(1 - {I_t \over n} \right) \left(1 - \sqrt{\eps}\right) \geq {\eps \over 2}
\left(1 - \sqrt{\eps}\right)  \geq {\eps \over 4}.$$ 
By substituting this bound into (\ref{eq:Recursion}) we obtain (\ref{eq:RecursioToProve}).

~\\
What remains is to show the second statement of the lemma. This follows readily from (\ref{eq:ILow}). Indeed, if $U_t$ denotes the number of uninformed 
vertices after $t$ stages, then observe first that $n\left(1 - {I_t / n} \right) = U_t$. So, for $n$ large enough
$$ U_{T_2(G)} = U_{T_2(G)-1} - N_{T_2(G)-1} \stackrel{(\ref{eq:ILow})}{\geq} 
U_{T_2(G)-1} e^{-I_{(T_2(G)}-1)/n}(1-e\sqrt{\eps})\geq {\eps n \over 2e}. $$
\end{proof}
Finally, we proceed by bounding $T(G) - T_2 (G)$, for $G \in \T_n (p)$. Let us denote by $\I_t$ the set of informed vertices after stage $t$. Recall that the main strategy in the previous argument was to show that the number $N_t$ of vertices that become informed by $\I_t$ in stage $t+1$ is close to its expected value. To achieve this, we exploited the fact that in $G$, except of a set $\X$ of size $\le \frac{8n}{\ln n}$, all vertices have the ``right'' degree in $\I_t$. This argument is unfortunately not applicable in the proof of the next lemma: for $t > T_2(G)$, the set $V \setminus \I_t$ of not yet informed vertices can become much smaller than $\X$, which makes our bounds useless. So we need to argue somehow differently.
An additional difficulty is that we are not able to apply the Azuma-Hoeffding inequality in a meaningful way. Note that for $t>T_2(G)$ the quantity $I_t$ 
is already of linear order, but the number $N_t$ of newly informed vertices at stage $t+1$ may become very small. In this case, the  Azuma-Hoeffding inequality gives a trivial upper bound and thus the need for a stronger concentration inequality. 
\begin{lem}
\label{lem:finishing}
Uniformly for all $G \in \T_n (p)$,  with probability $1-o(1)$
$$|(T(G) - T_2(G)) - \ln n| \le \eps^{1/3} \ln n.$$
\end{lem}
\begin{proof}
We will split the interval between $T_2(G)$ and $T(G)$ into two subintervals. In particular, let $T'(G)$ be the first time after $T_2(G)$ where at most $\ln^{1/2} n$ uninformed vertices remain. We will give separate bounds for $T'(G) - T_2 (G)$ and $T(G) - T'(G)$. Let for the remainder $\I_t$ be the random set of informed vertices after $t$ stages, and set $I_t := |\I_t|$. 

Let us start with the latter case, as it is the easier among the two. Let $t \ge T'(G)$. Since $np > \alpha (n) \ln n$, it follows from property (III)
that for every $v \in V \setminus \I_t$ we have for $n$ large enough
$|\Gamma_G (v) \cap \I_t| \geq np(1-3 \eps) - \ln^{1/2} n \geq  np(1-4\eps) $. 
So, the probability that a given uninformed vertex remains uninformed in the next stage is for large $n$ at most 
$$ \left( 1 - {1\over np(1+3\eps)}\right)^{np(1-4\eps)} \leq e^{-{1-4\eps \over 1+3\eps}} \leq {2\over e}.$$
Therefore, the probability that such a vertex remains uninformed for at least $\ln^{1/2} n$ 
steps after $T'(G)$ is at most $(2/e)^{\ln^{1/2} n}$. This implies that the expected number of vertices that remain uninformed for at least $\ln^{1/2} n$ stages after $T'(G)$ is at most $\ln^{1/2} n  (2/e)^{\ln^{1/2} n} \leq  
(2/e)^{\ln^{1/3} n} = o(1)$. That is, with probability at least $1 - (2/e)^{\ln^{1/3} n}$, we have $T(G) - T'(G) < \ln^{1/2} n$.

~\\
The bound on $T'(G) - T_2 (G)$ is significantly more complex. 
Let $U_t$ denote the number of vertices that are still uninformed after the $t$th stage. We will show that if $t$ is such that 
$U_t > \ln^{1/2} n$, then  
\begin{equation} \label{eq:ThirdRec}
U_{t+1} = U_t e^{-1} \left( 1\pm 50 \sqrt{\eps} \right),
\end{equation} 
with probability at least $1- e^{-\eps \ln^{1/2} n/10}$. 
So if $T'(G) - T_2(G) > \lceil \ln n + 55 \sqrt{\eps} \ln n \rceil =:b_1$, then with conditional probability at least
$(1- e^{-\eps \ln^{1/2} n/10})^{b_1+1} = 1-o(1)$ we have 
$$
U_{T_2(G) + b_1 +1}
\leq U_{T_2(G)} e^{-b_1-1} \left(1+ 50 \sqrt{\eps} \right)^{b_1+1}.$$
For large $n$
$$ \left(1+ 50 \sqrt{\eps} \right)^{b_1 + 1}\leq e^{55 \sqrt{\eps} \ln n}. $$
Also, $U_{T_2(G)} \leq \eps n$, which together with the above facts implies that $ U_{T_2(G) + b_1 +1} \le  \eps$. So, we may conclude that $T'(G) < T_2(G) + b_1 +1$. 

Similarly, if we assume that $T'(G) - T_2(G) < \lfloor \ln n - 55 \sqrt{\eps} \ln n \rfloor =:b_2$, then with conditional probability at least 
$(1- e^{-\eps \ln^{1/2} n/10})^{b_2} =1-o(1)$ we have
$$
U_{T_2(G) + b_2} \geq U_{T_2(G)} e^{-b_2} \left(1- 50 \sqrt{\eps} \right)^{b_2}.$$
A similar calculation as above, and the fact $U_{T_2(G)}\geq {\eps n \over 2e}$, which is guaranteed by Lemma~\ref{lem:slowDown} to hold with probability $1-o(1)$, shows that
$$U_{T_2(G) + b_2}\geq \eps e^{\sqrt{\eps} \ln n} \gg \ln^{1/2} n.$$
Thus,
$|T'(G) - T_2 (G) - \ln n| \leq   55 \sqrt{\eps} \ln n + 2$,
which concludes the proof of the lemma.

It remains to show~\eqref{eq:ThirdRec}.  
As an auxiliary preparation we will show that ``most'' vertices in $V \setminus \I_t$ have the ``right'' degree in $\I_t$, by arguing that if this was not the case, then there would be a significant deviation in the number of edges between $\I_t$ and $V \setminus \I_t$. More precisely, let
\[
	\X = \left\{v \in V\setminus \I_t : |\Gamma_G (v) \cap \I_t| < (1- 3 \sqrt{\eps})pn \right\}.
\]
In the sequel we argue that
\begin{equation}
\label{eq:XisSmall}
	|\X| \le 3\sqrt{\eps}(n - I_t).
\end{equation}
Indeed, as we assumed that $G\in\T_n(p)$, property (III) guarantees that $e_G(\I_t, V \setminus \I_t) \geq I_t (n- I_t)p (1- 3\eps)$.  
Moreover, property (III) implies that every vertex $v$ has degree at most $(1+3\eps)pn$. 
Therefore 
$$e_G(\I_t, V \setminus \I_t) < |\X|(1-3\sqrt{\eps})pn + (1+3\eps)(n-I_t-|\X|)pn .$$
By putting the upper and the lower bounds together we obtain
$$
I_t (n- I_t)p (1- 3\eps)
\leq
-3|\X|pn(\sqrt\eps + \eps) + (1+3\eps)(n-I_t)pn,
$$
which implies with $I_t \ge (1-\eps) n$ that 
$$1- 4\eps \leq  -3\frac{|\X|}{n - I_t}(\sqrt{\eps} + \eps) + (1+3\eps).$$
An elementary calculation shows that the claim \eqref{eq:XisSmall} holds.   

Now let $v \in (V\setminus \I_t) \setminus \X$. The probability that $v$ becomes informed in the next stage is 
$$1 - \prod_{u \in \Gamma_G (v) \cap \I_t} \left(1 - {1\over |\Gamma_G (u)|} \right) = 
1 - \left( 1 - {1 \over pn(1\pm 3\eps)}\right)^{pn (1\pm 3\sqrt{\eps})}   
= 1 - {1\over e} \left(1\pm 7\sqrt{\eps} \right).
$$ 
Denote by $\N_t$ the set of vertices in $V\setminus \I_t$ that will be informed by the vertices in $\I_t$ in stage $t+1$. Moreover, write $N_t = |\N_t|$. So, by linearity of expectation, for $n$ large enough we obtain
\begin{equation} \label{eq:ExpLowThird}
\begin{split}
\ex \left( \RN \right) & =  
(n - I_t)  \left(1 - {1\over e} \right)(1-7\sqrt{\eps}) \pm  3\sqrt{\eps} (n- I_t) \\
& = (n - I_t)  \left(1 - {1\over e} \right)(1\pm  14\sqrt{\eps}).
\end{split}
\end{equation}
Next we will show that $\RN$ is with sufficiently high probability close to its expected value. Note that the Azuma-Hoeffding inequality does not give any meaningful bounds, as the number of the independent random variables is $I_t \ge (1-\eps)n$, while the expected value of $\RN$ is proportional only to $n- I_t$. The latter 
will eventually become so small that the exponent in the Azuma-Hoeffding inequality is $o(1)$, thus yielding a trivial bound.  To bypass this problem, we will use Talagrand's inequality (Theorem~\ref{thm:talagrand}). Note first that the bounded differences condition is satisfied, that is, changing one random choice can change $\RN$ by at most 1. Regarding the second condition, note that if $\RN=r$, then there must be at  least $r$ vertices in $\I_t$ that have informed the vertices in $\N_t$. Therefore, we may take $\psi (r) = \lceil r\rceil$ and with $m(\RN)$ denoting the median 
of $\RN$, we deduce for any $x>0$ that 
\begin{equation} \label{eq:TalUltim}
\prob \left(|\RN - m(\RN)| > x \right)\leq 4\exp \left(- {x^2 \over 4(m(\RN) + x)} \right) \leq 
4\exp \left(- {x^2 \over 4(2\ex(\RN) + x)} \right),
\end{equation}
where in the last inequality we have used the fact that $\ex (\RN) \geq m (\RN ) \prob (\RN >  m (\RN))\geq m (\RN )/2$, 
which implies that $m(\RN) \leq 2 \ex (\RN)$. 
However, we need to argue about the distance of $m(\RN)$ from $\ex (\RN)$. 
We will use (\ref{eq:MedMean}).
The triangle inequality yields:
\begin{equation*} 
\begin{split}
& | \RN - \ex (\RN)| = |\RN - m (\RN) + m(\RN) - \ex (\RN) | \\
& \leq |\RN - m (\RN) | + |\RN - m (\RN) | \stackrel{(\ref{eq:MedMean})}{=} 
|\RN - m (\RN) |+ O \left(\sqrt{\ex (\RN )} \right).
\end{split}
\end{equation*} 
Since $\alpha (n) \leq \ln^{1/9}n$, we have $\sqrt{\ex (\RN)} = o(\sqrt{\eps} \ex (\RN))$. Therefore, for sufficiently large $n$
\begin{equation} \label{eq:NewTail} 
|\RN - \E{\RN}) | > x
\implies
|\RN - m (\RN) | > x - \sqrt{\eps} \ex (\RN).
\end{equation} 
Therefore, using (\ref{eq:NewTail}) in (\ref{eq:TalUltim}) with $x = \sqrt{\eps}\E{N_t}$ we obtain
\begin{equation} \label{eq:TalAppl}
\begin{split}
\prob ( |\RN - \ex ( \RN)| > 2\sqrt{\eps} \ex (\RN)) \leq 4 \exp \left(- {\eps \ex (\RN) \over 4(2+\sqrt{\eps})} \right). 
\end{split}
\end{equation}
Since $n- I_t \geq \ln^{1/2} n$, by (\ref{eq:ExpLowThird}) we obtain that, say, 
$\ex (\RN) \geq {\ln^{1/2} n \over 3}$. So, for large $n$, the bound in (\ref{eq:TalAppl}) becomes 
\begin{equation*} 
\begin{split}
\prob ( |\RN - \ex ( \RN)| > 2\sqrt{\eps} \ex (\RN)) \leq \exp \left(- {\eps \ln^{1/2} n \over 40} \right). 
\end{split}
\end{equation*}
By putting everything together we obtain that with probability at least $1 - e^{-\eps \ln^{1/2} n /40}$ 
$$ \RN = (n - I_t)  \left(1 - {1\over e} \right) \left( 1\pm 16\sqrt{\eps}\right). $$
So there remain (very generously)
$ (n-I_t) e^{-1} \left(1\pm 60 \sqrt{\eps}  \right)$
vertices uninformed in $V \setminus \I_t$. This completes the proof of \eqref{eq:ThirdRec}.
\end{proof}
Finally, note that the bounds obtained in Lemmas~\ref{lem:doubling}--\ref{lem:finishing} imply Theorem~\ref{thm:main}, thus concluding 
our proof.

\bibliographystyle{plain}
\bibliography{RandBroad}

\end{document}